# FINITE FOURIER SERIES. THE CLASS OF TRIGONOMETRIC SPLINES

**Denysiuk V.P.**[1], Doctor of Physical and Mathematical Sciences, Professor

**Rybachuk L.V.**[2], Candidate of Physical and Mathematical Sciences (Ph.D.), Docent

[1] National Aviation University, Ukraine

[2] National Technical University of Ukraine "Igor Sikorsky Kyiv Polytechnic Institute", Ukraine

(kvomden@ukr.net, rybachuk.liudmyla@lll.kpi.ua)

### Abstract

Finite trigonometric Fourier series on a set of discrete equidistant points are considered. A finite system of orthogonal functions that have interpolation and certain differential properties on the period $[0, 2\pi)$ is introduced. Finite Fourier series based on this system of functions form a class of trigonometric splines, which includes polynomial periodic simple splines. Trigonometric splines suggest generalizations in several directions and certainly require further research.

### Introduction

Trigonometric Fourier series have played an important role in the development of mathematics. Along with the conventional theory, in which the function $f(t)$ is considered periodic with period $2\pi$ and known over the entire period, there is a theory that considers only discrete values of a certain function on a set of equidistant points; in this case, the expansion is based on a finite system of functions of the form

$$1, \cos t, \sin t, \cos 2t, \sin 2t, \ldots, \cos nt, \sin nt . \qquad (1)$$

Such series are called finite Fourier series [1].

The sums of such series can be considered as analytical functions that interpolate a given function $f(t)$ at the nodes of uniform grids. However, in many cases, it is necessary to construct an interpolation function with specified differential properties. It is clear that such a task requires the construction of another finite system of functions that coincides with (1) at the nodes of uniform grids and has the specified differential properties.

### Purpose of the work

Construction of a finite system of functions that coincides with system (1) at the nodes of uniform grids and has specified differential properties.

### The main part

Let us consider the theory of finite Fourier series in more detail, following [1].

Assume that on the interval $[0, 2\pi)$ a uniform grid $\Delta_N = \{t_j\}_{j=1}^{N}$, $t_j = \frac{2\pi}{N}(j-1)$, $N = 2n+1$, $n = 1, 2, \ldots$ is given.

The expansion in a Fourier series is based on the orthogonality of the functions

$$1, \cos t, \sin t, \cos 2t, \sin 2t, \ldots, \cos kt, \sin kt, \ldots$$

with respect to the integration operation on the interval $[0, 2\pi]$.

$$\int_0^{2\pi} \cos kt \cos mt\, dt = \begin{cases} 0, & k \neq m; \\ \pi, & k = m \neq 0; \\ 2\pi, & k = m = 0; \end{cases}$$

$$\int_0^{2\pi} \sin kt \sin mt\, dt = \begin{cases} 0, & k \neq m; \\ \pi, & k = m \neq 0; \end{cases}$$

$$\int_0^{2\pi} \sin kt \cos mt\, dt = 0.$$

If, instead of the integration operation, the summation operation is applied, then the functions (1) turn out to be orthogonal on the discrete set of points that are nodes of the grids $\Delta_N$. Indeed, it is easy to verify that the following relations hold

$$\sum_{j=1}^{N} \cos kt_j \cos mt_j = \begin{cases} 0, & k \neq m; \\ N/2, & k = m \neq 0; \\ N, & k = m = 0; \end{cases} \qquad (2)$$

$$\sum_{j=1}^{N} \sin kt_j \sin mt_j = \begin{cases} 0, & k \neq m; \\ N/2, & k = m \neq 0; \end{cases}$$

$$\sum_{j=1}^{N} \sin kt_j \cos mt_j = 0.$$

These relations can be considered as analogs of orthogonality conditions [1].

Let us consider a trigonometric polynomial of the order $n$ ($n = 1, 2, \ldots$) of the form

$$T_n(t) = \frac{a_0}{2} + \sum_{k=1}^{n} a_k \cos kt + b_k \sin kt. \qquad (3)$$

The right-hand side of (3) is called the finite Fourier series, and the polynomial $T_n(t)$ is the sum of this series [1]. Note that the polynomial $T_n(t)$ is considered on a discrete set of points $\Delta_N$.

Let us represent the values of some discrete function $\{f(t_j)\}_{j=1}^{N} = \{f_j\}_{j=1}^{N}$ at the nodes of the grid $\Delta_N$ by a trigonometric polynomial (3):

$$f_j = T_n(t_j) = \frac{a_0}{2} + \sum_{k=1}^{n} a_k \cos kt_j + b_k \sin kt_j, \qquad (4)$$

where $a_0, a_k, b_k$ ($k = 1, 2, \ldots, n$) are some coefficients. To determine these coefficients, we will use the orthogonality relation (2).

The coefficient $a_0$ is determined as follows. Multiply (3) by $\cos 0$ and sum both parts over all $j$ ($j = 1, 2, \ldots, N$). We have

$$\sum_{j=1}^{N} f_j^{(0)} \cos 0 = \frac{a_0}{2} \sum_{j=1}^{N} \cos 0 + \sum_{k=1}^{n} a_k \sum_{j=1}^{N} \cos 0 \cos kt_j + \sum_{k=1}^{n} b_k \sum_{j=1}^{N} \cos 0 \sin kt_j.$$

Taking into account (2), we obtain:

$$\sum_{j=1}^{N} f_j \cos 0 = \frac{a_0}{2} N.$$

Hence,

$$a_0^{(0)} = \frac{2}{N} \sum_{j=1}^{N} f_j . \tag{5}$$

The coefficients of $a_k$ and $b_k$ ($k = 1, 2, ..., n$) are determined in the same way. Multiply (3) by $\cos kt_j$ and by $\sin kt_j$, and sum both sides over all $j$ ($j = 1, 2, ..., N$). Taking into account (2), we obtain, respectively:

$$\sum_{j=1}^{N} f_j \cos kt_j = a_k \frac{N}{2} ; \quad \sum_{j=1}^{N} f_j \sin kt_j = b_k \frac{N}{2},$$

from which

$$a_k = \frac{2}{N} \sum_{j=1}^{N} f_j \cos kt_j ; \quad b_k = \frac{2}{N} \sum_{j=1}^{N} f_j \sin kt_j , \quad k = 1, 2, ..., n . \tag{6}$$

Let us consider an example. Let $N = 9$ and $\{f_j\}_{j=1}^{N} = \{2, 1, 3, 2, 4, 1, 3, 1, 3\}$. The graph of the discrete trigonometric polynomial $T_4(t)$ is shown in Figure 1. Note that here and in the following, the vertical lines on the graphs coincide with the nodes of the grid $\Delta_N$.

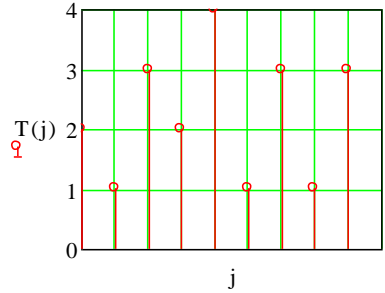

Fig. 1. Discrete trigonometric polynomial $T_4(t)$ on the grid $\Delta_N$

Until now, we have considered the discrete function $\{f_j\}_{j=1}^{N}$ and the trigonometric polynomial only at the nodes of the grid $\Delta_N$. In the future, we will consider the continuous function $f(t)$, which takes the value $\{f_j\}_{j=1}^{N}$ at the nodes of the grid $\Delta_N$, and the trigonometric polynomial (3) on the interval $[0, 2\pi)$. In this case, the coefficients (5), (6) can be obtained based on other considerations (see, e.g., [2]).

Let us consider an example. Let $N = 9$ and $\{f_j\}_{j=1}^{N} = \{2, 1, 3, 2, 4, 1, 3, 1, 3\}$ be the same as before. The graph of the polynomial $T_4(t)$ is shown in Figure 2.

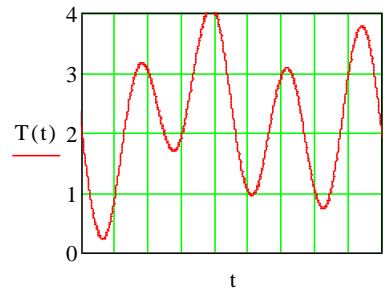

Fig. 2. Trigonometric polynomial $T_4(t)$ with coefficients (5), (6),
which interpolates the function $f(t)$ at the nodes of the grid $\Delta_N$

It is clear that the polynomial $T_n(t)$ with coefficients (4), (5) is an analytic function. Naturally, the question arises about the possibility of constructing a system of functions that satisfy the conditions:

a) at the nodes of the grid $\Delta_N$, these functions coincide with the system (1);

b) on the interval $[0, 2\pi)$, these functions have the given differential properties.

It turns out that such systems of functions can be constructed, for example, in the following way.

Consider the functions

$$C_k(\sigma, r, N, t) = \frac{C(\sigma, r, N, k, t)}{H(\sigma, r, N, k)}, \quad S_k(\sigma, r, N, t) = \frac{S(\sigma, r, N, k, t)}{H(\sigma, r, N, k)} \quad (k = 1, 2, ..., n), \tag{7}$$

where

$$\sigma(r, k) = k^{-(1+r)}; \tag{8}$$

$$C(\sigma, r, N, k, t) = \sigma(r, k) \cos kt +$$
$$+ \sum_{m=1}^{\infty} \left[ \sigma(r, mN + k) \cos((mN + k)t) + (-1)^{1+r} \sigma(r, mN - k) \cos((2mN - k)t) \right]; \tag{9}$$

$$S(\sigma, r, N, k, t) = \sigma(r, k) \sin kt +$$
$$+ \sum_{m=1}^{\infty} \left[ \sigma(r, mN + k) \sin((mN + k)t) - (-1)^{1+r} \sigma(r, mN - k) \sin((mN - k)t) \right]; \tag{10}$$

$$H(\sigma, r, N, k) = \sigma(r, k) + \sum_{m=1}^{\infty} \left[ \sigma(r, mN + k) + (-1)^{1+r} \sigma(r, mN - k) \right]; \tag{11}$$

$$(k = 1, 2, ..., n; \quad r = 1, 2, ...).$$

It is easy to verify that for any value of the parameter $r$, the following equalities hold

$$C_k(\sigma, r, N, t_j) = \cos kt_j; \quad S_k(\sigma, r, N, t_j) = \sin kt_j \quad (j = 1, 2, ... N).$$

Indeed,

$$\cos(mN + k)t_j = \cos(mNt_j + kt_j) = \cos\left[mN \frac{2\pi}{N}(j-1) + k \frac{2\pi}{N}(j-1)\right] =$$
$$= \cos\left[2\pi m(j-1) + k \frac{2\pi}{N}(j-1)\right] = \cos k \frac{2\pi}{N}(j-1) = \cos kt_j.$$

$$\cos(mN - k)t_j = \cos kt_j.$$

Similarly

$$\sin(mN + k)t_j = \sin kt_j; \quad \sin(mN - k)t_j = -\sin kt_j.$$

Then

$$C(\sigma, r, N, k, t_j) = \sigma(r, k) \cos kt_j +$$
$$+ \sum_{m=1}^{\infty} \left[ \sigma(r, mN + k) \cos((mN + k)t_j) + (-1)^{1+r} \sigma(r, mN - k) \cos((2mN - k)t_j) \right] =$$
$$= \cos kt_j \left[ \sigma(r, k) + \sum_{m=1}^{\infty} \left[ \sigma(r, mN + k) + (-1)^{1+r} \sigma(r, mN - k) \right] \right] = \cos kt_j H(\sigma, r, k).$$

$$S(\sigma, r, N, k, t_j) = \sigma(r, k)\sin kt_j +$$

$$+ \sum_{m=1}^{\infty} \left[ \sigma(r, mN+k)\sin((mN+k)t_j) - (-1)^{1+r}\sigma(r, mN-k)\sin((2mN-k)t_j) \right] =$$

$$= \sin kt_j \left[ \sigma(r, k) + \sum_{m=1}^{\infty} \left[ \sigma(r, mN+k) + (-1)^{1+r}\sigma(r, mN-k) \right] \right] = \sin kt_j H(\sigma, r, k).$$

Finally, we have

$$C_k(\sigma, r, N, t_j) = \frac{C(\sigma, r, Nk, t_j)}{H(\sigma, r, N, k)} = \frac{\cos t_j H(\sigma, r, N, k)}{H(\sigma, r, N, k)} = \cos kt_j \,;$$

$$S_k(\sigma, r, N, t_j) = \frac{S(\sigma, r, N, k, t_j)}{H(\sigma, r, N, k)} = \frac{\sin t_j H(\sigma, r, N, k)}{H(\sigma, r, N, k)} = \sin kt_j \quad (j = 1, 2, ..., N).$$

Thus, the functions $C_k(\sigma, r, N, t)$, $S_k(\sigma, r, N, t)$ ($j = 1, 2, ..., n$), take the same values at the nodes of the grid $\Delta_N$ as the system (1) for any values of the parameter $r$ ($r = 1, 2, ...$).

Let us now consider the differential properties of the functions $C(\sigma, r, N, k, t)$ and $S(\sigma, r, N, k, t)$. First of all, note that the sign-constant factors $\sigma(r, k)$ have the order of decay $O(k^{-(1+r)})$. It is clear that the infinite trigonometric series representing these functions converge uniformly at $r = 1$ (according to the Weierstrass criterion), and, therefore, their sums are continuous functions. For $r > 1$, these series and their derivatives of order $r - 1$ also converge uniformly, and, therefore, belong to the class $C_p^{r-1}$, where $C_p^{r-1}$ is the class of $2\pi$-periodic functions that have continuous derivatives up to order $r - 1$ inclusive. Since constant factors do not affect the differential properties, we can conclude that the functions $C_k(\sigma, r, N, t)$ and $S_k(\sigma, r, N, t)$ also belong to the class $C_p^{r-1}$.

Thus, we have constructed the trigonometric functions $C_k(\sigma, r, N, t)$ and $S_k(\sigma, r, N, t)$, which take the same values at the nodes of the grids as system (1) and have certain differential properties — they belong to the class $C_p^{r-1}$. Therefore, from the system of functions (1) we can move to the system

$$1, C_1(\sigma, r, N, t), S_1(\sigma, r, N, t), ..., C_n(\sigma, r, N, t), S_n(\sigma, r, N, t). \tag{12}$$

Let us consider in more detail the functions of type (3) with the components $C_k(\sigma, r, N, t)$ and $S_k(\sigma, r, N, t)$ of the form

$$T_n^*(\sigma, r, f, N, t) = \frac{a_0}{2} + \sum_{k=1}^{n} a_k C_k(\sigma, r, N, t) + b_k S_k(\sigma, r, N, t). \tag{13}$$

It is clear that the functions $T_n^*(\sigma, r, f, N, t)$ belong to the class $C_p^{r-1}$ and interpolate the function $f(t)$, the values of which were used to construct the coefficients $a_0, a_k, b_k$ ($k = 1, 2, ..., n$) at the nodes of the grid $\Delta_N$. Since periodic polynomial splines have the same properties (see, for example, [6]), it is natural to call the functions $T_n^*(\sigma, r, f, N, t)$ trigonometric splines; these splines will be denoted by $St(\sigma, r, f, N, t)$. Hence,

$$St(\sigma, r, f, N, t) = \frac{a_0}{2} + \sum_{k=1}^{n} a_k C_k(\sigma, r, t) + b_k S_k(\sigma, r, t). \tag{14}$$

For illustration, we will give some examples of splines $St(\sigma, r, f, N, t)$ for different values of the parameter $r$ ($r = 0, 1, ...$); also for comparison, we will give a graph of the polynomial $T_4(t)$. Note that although the 0th order splines are not splines in the sense of the definitions introduced above (the conditions for uniform convergence of the trigonometric series are not satisfied), their consideration is still very useful in many cases.

Let, as before, $N = 9$ and $\{f_j\}_{j=1}^{N} = \{2,1,3,2,4,1,3,1,3\}$. Note that on the graphs below, the splines of $St(\sigma, r, f, N, t)$ are denoted briefly, namely $St(r, t)$ (Figure 3).

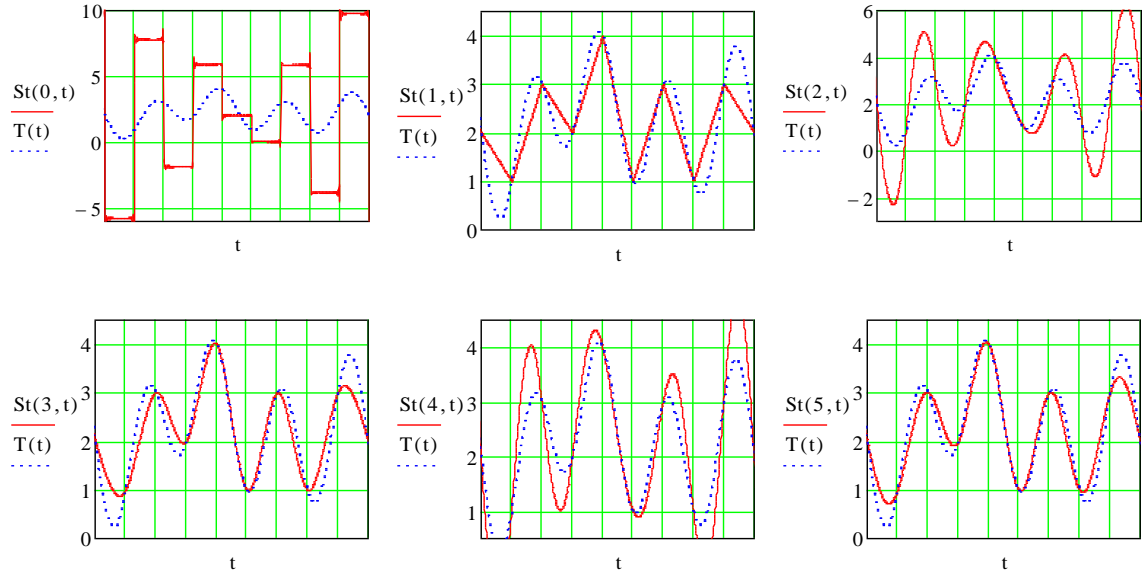

Fig. 3. Graphs of trigonometric splines $St(\sigma, r, f, N, t)$; $r = 0, 1, ..., 5$

The presented trigonometric splines allow for broad generalizations, which we will consider in more detail.

Earlier, on the interval $[0, 2\pi)$ we considered the equidistant grid $\Delta_N = \{t_j\}_{j=1}^{N}$, where $t_j = \frac{2\pi}{N}(j-1)$, $N = 2n+1$, $n = 1, 2, ...$. Further, this grid will be denoted as $\Delta_N^{(0)}$ and its nodes as $\{t_j^{(0)}\}_{j=1}^{N}$. We will also consider the grid $\Delta_N^{(1)} = \{t_j^{(1)}\}_{j=1}^{N}$, $t_j^{(1)} = \frac{\pi}{2N}(2j-1)$, ($N = 2n+1$, $n = 1, 2, ...$), given on $(0, 2\pi)$. Let us also introduce the indicator $I$ ($I = 0, 1$); then both grids can be denoted as $\Delta_N^{(I)}$, and their nodes as $\{t_j^{(I)}\}_{j=1}^{N}$.

It is easy to verify that the relations of type (2) are also fulfilled at the nodes of the grid $\Delta_N^{(1)}$. Then formulas (4), (5), (6) can be written in the form

$$f_j^{(I)} = T_n^{(I)}(t_j^{(I)}) = \frac{a_0^{(I)}}{2} + \sum_{k=1}^{n} a_k^{(I)} \cos k t_j^{(I)} + b_k^{(I)} \sin k t_j^{(I)}, \qquad (15)$$

where

$$a_0^{(I)} = \frac{2}{N} \sum_{j=1}^{N} f_j^{(I)};$$

$$a_k^{(I)} = \frac{2}{N} \sum_{j=1}^{N} f_j^{(I)} \cos k t_j^{(I)}; \quad b_k^{(I)} = \frac{2}{N} \sum_{j=1}^{N} f_j^{(I)} \sin k t_j^{(I)}, \quad k = 1, 2, ..., n. \qquad (16)$$

Let us introduce the concept of stitching and interpolation grids.

The stitching grid of trigonometric splines $St(\sigma, r, t)$ will be called the grid at the nodes of which discontinuities of the function itself occur (for splines of order 0) or discontinuities of the derivatives of order $r$ (for splines of order $r = 1, 2, ...$).

The interpolation grid of trigonometric splines will be called the grid at the nodes of which the interpolation of a given function occurs.

Stitching and interpolation grids are best distinguished when considering the 0th and 1st order splines.

Thus, for example, from Figure 3, it can be concluded that in the presented splines, the stitching and interpolation grids coincide regardless of the evenness of their degree.

Let us introduce trigonometric splines that depend on several parameters:

$$St(I1, I2, I3(\sigma), \Lambda, \Gamma, H, \sigma, r, f, N, t) =$$

$$= \frac{a_0}{2} + \sum_{k=1}^{n} \lambda_k \left[ C_k(I1, I2, I3(\sigma), \Gamma, \sigma, r, t) a_k^{1-I2} a1_k^{I2} + S_k(I1, I2, I3(\sigma), H, \sigma, r, t) b_k^{1-I2} b1_k^{I2} \right], \quad (17)$$

where the parameters

$I1$ ( $I1 = 0,1$ ) determine the type of stitching grid;

$I2$ ( $I2 = 0,1$ ) determine the type of interpolation grid;

$I3(\sigma, r)$ ( $I3(\sigma, r) = 0,1$ ) determine the type of spline of even degree (the role of this parameter will be discussed in more detail below);

$\Lambda$ ( $\Lambda = (\lambda_1, \lambda_2, ..., \lambda_n)$ ) — are the chosen factors of the $\Lambda$-summation method;

$\Gamma$ ( $\Gamma = (\gamma_1, \gamma_2, \gamma_3)$., $\gamma_2 \gamma_3 \neq 0$ ) and $H$ ( $H = (\eta_1, \eta_2, \eta_3)$, $\eta_2 \eta_3 \neq 0$ ) are vectors of parameters that determine the influence of the low-frequency, mid-frequency, and high-frequency components of the spline;

$r$ ( $r = 1, 2, ...$ ) is the order of the trigonometric spline;

$\sigma(r, k)$ is the type of convergence factor which has an order of decay $O(k^{-(1+r)})$;

$$C_k(I1, I2, I3(\sigma), \Gamma, \sigma, r, t) = \gamma_1 \sigma(r, k) \cos kt +$$

$$+ \sum_{m=1}^{\infty} (-1)^{mK1(\sigma, r)} \left[ \gamma_3 \sigma(r, mN + k) \cos((mN + k)t) + (-1)^{(1+r)I3(\sigma)} \gamma_2 \sigma(r, mN - k) \cos((mN - k)t) \right]; \quad (18)$$

$$S_k(I1, I2, I3(\sigma), H, \sigma, r, t) = \eta_1 \sigma(r, k) \sin kt +$$

$$+ \sum_{m=1}^{\infty} (-1)^{mK1(\sigma, r)} \left[ \eta_3 \sigma(r, mN + k) \sin((mN + k)t) - (-1)^{(1+r)I3(\sigma)} \eta_2 \sigma(r, mN - k) \sin((mN - k)t) \right]; \quad (18a)$$

$$Hc_k(I1, I2, I3(\sigma), \Gamma, \sigma, r) = \gamma_1 \sigma(r, k) +$$

$$+ \sum_{m=1}^{\infty} (-1)^{mK2(\sigma, r)} \left[ \gamma_3 \sigma(r, mN + k) + (-1)^{(1+r)I3(\sigma)} \gamma_2 \sigma(r, mN - k) \right]; \quad (19)$$

$$Hs_k(I1, I2, I3(\sigma), H, \sigma, r) = \eta_1 \sigma(r, k) +$$

$$+ \sum_{m=1}^{\infty} (-1)^{mK2(\sigma, r)} \left[ \eta_3 \sigma(r, mN + k) + (-1)^{(1+r)I3(\sigma)} \eta_2 \sigma(r, mN - k) \right]. \quad (19a)$$

In the following, the set of trigonometric splines represented by expression (15) will be called the class of trigonometric splines and denoted by the expression $ST(f, r, N)$.

The coefficients $K1(\sigma, r)$, $K2(\sigma, r)$ depend on the chosen type of convergence factors $\sigma(r, k)$. In [3], the case was considered when the sign-changing Riemann factors [4] were chosen as convergence factors, which have the form

$$\sigma(r, k) = \left( \frac{\sin(\frac{\pi}{N} k)}{\frac{\pi}{N} k} \right)^{1+r} \quad (k = 1, 2, ...). \quad (20)$$

and in [5], sign-constant convergence factors were considered

$$\sigma(r,k) = k^{-(1+r)} \qquad (k = 1, 2, \ldots). \tag{21}$$

In the following, the Riemann convergence factors (20) will be denoted as $\sigma 0(r,k)$, and the convergence factors of the form (21) will be denoted as before by $\sigma(r,k)$. Accordingly, for the sign-changing convergence factor $\sigma 0(r,k)$, the coefficients $K1(\sigma 0, r)$, $K2(\sigma 0, r)$ are as follows

$$K1(\sigma 0, r) = r + 1 + I1, \quad K2(\sigma 0, r) = r + 1 + I1 + I2, \tag{22}$$

and for the sign-constant convergence factors $\sigma(r,k)$, the coefficients $K1(\sigma, r)$, $K2(\sigma, r)$ are of the form

$$K1(\sigma, r) = I1, \quad K2(\sigma, r) = I1 + I2. \tag{23}$$

Thus, by choosing the convergence factors $\sigma 0(r,k)$ or $\sigma(r,k)$, we also set the corresponding expressions for the coefficients $K1(\sigma 0, r)$, $K2(\sigma 0, r)$ and $K1(\sigma, r)$, $K2(\sigma, r)$. Therefore, in the notation of trigonometric splines, we do not indicate the dependence on these coefficients.

The question of the relation between the classes of polynomial periodic interpolation splines, whose theory is well developed (see, e.g., [6], [7], [8]) and the introduced class of trigonometric interpolation splines $ST(f, r, N)$, is very interesting. Let us consider this issue in more detail.

First of all, note that the trigonometric splines with Riemann convergence factors (20) and $I3(\sigma 0) = 0$, coincide with the trigonometric splines with convergence factors (19) and $I3(\sigma) = 1$, i.e.

$$St(I1, I2, 0, \Lambda, \Gamma, H, \sigma 0, r, f, N, t) \equiv St(I1, I2, 1, \Lambda, \Gamma, H, \sigma, r, f, N, t). \tag{24}$$

Similarly,

$$St(I1, I2, 1, \Lambda, \Gamma, H, \sigma 0, r, f, N, t) \equiv St(I1, I2, 0, \Lambda, \Gamma, H, \sigma, r, f, N, t). \tag{25}$$

Let us consider the unit vectors $\Lambda 1$ ($\Lambda 1 = (1,1,\ldots,1)$), $\Gamma 1$ ($\Gamma 1 = (1,1,1)$), and $H1$ ($H1 = (1,1,1)$). Trigonometric splines $St(I1, I2, I3(\sigma), \Lambda 1, \Gamma 1, H1, \sigma, r, f, N, t)$ with unit vectors will be denoted as $St(I1, I2, I3(\sigma), \sigma, r, f, N, t)$.

In [9] it was shown that trigonometric splines $St(0, 0, 0, \sigma 0, 2k+1, f, N, t)$ ($k = 0, 1, \ldots$) of odd orders for all $k$ coincide with periodic polynomial splines of the same powers $2k + 1$. Trigonometric splines $St(0, 0, 1, \sigma 0, 2k+1, f, N, t)$ ($k = 0, 1, \ldots$) of odd orders for all $k$ also coincide with periodic polynomial splines of the same degree $2k + 1$.

Splines of even order $St(0, 1, 0, \sigma 0, 2(k+1), f, N, t)$ for all $k$ ($k = 0, 1, \ldots$) coincide with periodic polynomial splines with a stitching grid $\Delta_N^{(0)}$ and an interpolation grid $\Delta_N^{(1)}$.

Taking into account (24) and (25), similar conclusions can be transferred to the trigonometric splines $St(0, 0, 1, \sigma, r, f, N, t)$ and $St(0, 0, 0, \sigma, r, f, N, t)$ with a sign-constant convergence factor (21).

The presented conclusions are difficult to overestimate, as they allow transferring all the results obtained in the approximation theory for polynomial splines to the indicated types of trigonometric splines. In addition, the established connection between the classes of polynomial and trigonometric splines allows applying the methods of trigonometric series theory to the study of polynomial spline classes and vice versa.

The empirical rule by which the existence of polynomial analogs of trigonometric splines of arbitrary degree $r$ ($r = 1, 2, \ldots$) is established can be formulated as follows.

Rule. If the derivative of the $r$-th order of the trigonometric spline of the $r$-th order is a piecewise constant (step) function, then there exists a polynomial analog of such a spline. If such a derivative is not a piecewise constant function, then the trigonometric spline does not have a polynomial analog.

Thus, for example, splines of even order $St(0, 0, 1, \sigma 0, 2(k+1), f, N, t)$ ($k = 0, 1, \ldots$), which have stitching and interpolation grids $\Delta_N^{(0)}$, do not have polynomial analogues.

Polynomial analogs of trigonometric splines with other parameters are currently unknown, although in some cases they can be easily constructed (at least the first-order splines with unit vectors $\Lambda 1$, $\Gamma 1$ and $H1$, which are broken lines).

Trigonometric splines of the class $ST(f,r,N)$ allow various forms of representation. For example, splines can be represented by fundamental splines [10]; using them, it is easy to construct trigonometric splines of several variables [11]. It is also possible to represent trigonometric splines through trigonometric $B$-splines [12].

One of the advantages of trigonometric splines over polynomial splines is that trigonometric splines are represented by a single expression, while polynomial splines are composite functions; this significantly simplifies the performance of mathematical analysis operations on trigonometric splines. Another advantage of trigonometric splines is the significant simplification of constructing these splines of arbitrary orders.

Trigonometric splines allow for generalizations in several directions. For example, these splines can be used to approximate non-periodic functions [13]. Further, by choosing the parameters $I1$ and $I2$, various combinations of stitching and interpolation grids can be obtained. It is also possible to set non-unit vectors $\Gamma$ and $H$ in these splines, thereby changing the ratio between the low-frequency, mid-frequency, and high-frequency components of these functions. The possibility of applying various types of convergence factors, the only requirement for which is to ensure the order of decay $O(k^{-(1+r)})$ [14], also deserves attention. Finally, trigonometric splines with fractional derivatives in the sense of Weyl [15] can also be considered.

Undoubtedly, the class of trigonometric splines requires further research.

## Conclusions

1. A finite system of functions

    $$1, C_k(I1, I2, I3(\sigma), \Gamma, \sigma, r, t), S_k(I1, I2, I3(\sigma), H, \sigma, r, t) \quad (k = 1, 2, ..., n),$$

    which coincides with the system

    $$1, \cos t, \sin t, \cos 2t, \sin 2t, ..., \cos nt, \sin nt,$$

    at the nodes of the grids $\Delta_N^{(I)}$ and belongs to the class $C_p^{r-1}$ ($r = 1, 2, ...$), has been constructed.

2. Finite Fourier series based on this system of functions form the class of trigonometric splines $ST(f,r,N)$.

3. Trigonometric interpolation splines of the class $ST(f,r,N)$ in some cases coincide with periodic polynomial interpolation splines; this fact allows transfering all the results obtained in the theory of approximations for polynomial splines to the indicated types of trigonometric splines.

4. The established connection between the classes of polynomial and trigonometric splines allows applying the methods of trigonometric series theory to study the classes of polynomial splines and vice versa.

5. The trigonometric splines of the class $ST(f,r,N)$ allow for generalizations in several directions.

6. Trigonometric splines of the class $ST(f,r,N)$ can be represented by fundamental trigonometric splines and by trigonometric $B$-splines.

7. The construction of trigonometric splines of arbitrary orders is significantly simpler than the construction of corresponding polynomial splines of the same degree.

8. The class $ST(f,r,N)$ contains trigonometric splines that do not have polynomial analogs.

9. Trigonometric splines of the class $ST(f,r,N)$ allow for the possible application of non-integer values of the parameter $r$ ($r > 0$), at which splines and their derivatives of fractional orders are formed.

10. Undoubtedly, the class of trigonometric splines $ST(f,r,N)$ requires further research.